%
%

\magnification=1200


\font\AAA=cmr14 at 12pt
\font\BBB=cmr14 at 8pt

\overfullrule=0in

\def\boxit#1{\hbox{\vrule
 \vtop{%
  \vbox{\hrule\kern 2pt %
     \hbox{\kern 2pt #1\kern 2pt}}%
   \kern 2pt \hrule }%
  \vrule}}

  \def\Arr#1{\buildrel {#1} \over \longrightarrow}

\def\pn{\bbp^n}

\def\ss{\subset}

\def\smfrac#1#2{\hbox{${#1\over #2}$}}

\def\deg{{\rm deg}}

\def\log{{\rm log}}

\def\max{{\rm max}}
\def\min{{\rm min}}

\def\arr{\longrightarrow}
\def\supp{{\rm supp}}


\def\Theorem#1{\medskip\noindent {\AAA T\BBB HEOREM \rm #1.}}
\def\Prop#1{\medskip\noindent {\AAA P\BBB ROPOSITION \rm  #1.}}
\def\Cor#1{\medskip\noindent {\AAA C\BBB OROLLARY \rm #1.}}
\def\Lemma#1{\medskip\noindent {\AAA L\BBB EMMA \rm  #1.}}

\def\Note#1{\medskip\noindent {\AAA N\BBB OTE \rm  #1.}}
\def\Def#1{\medskip\noindent {\AAA D\BBB EFINITION \rm  #1.}}

\def\Conj#1{\medskip\noindent {\AAA C\BBB ONJECTURE \rm    #1.}}

\def\pf{\medskip\noindent {\bf Proof.}\ }
\def\qed{\hfill  $\vrule width5pt height5pt depth0pt$}

   \def\cp{{\cal P}}   
   \def\co{{\cal O}}

\def\cl{{\cal L}}
\def\cp{{\cal P}}

\def\wh{\widehat}

\def\and{\qquad {\rm and} \qquad}
\def\arr{\longrightarrow}

\def\bbc{{\bf C}}

\def\bbp{{\bf P}}

\def\a{\alpha}

\def\d{\delta}

\def\f{\phi}
\def\g{\gamma}

\def\l{\lambda}

\def\s{\sigma}

\def\z{\zeta}

\def\D{\Delta}
\def\L{\Lambda}

\def\O{\Omega}

\def\PH#1{\widehat {#1}}

\def\bo{\partial \Omega}

\def\AA{1}
\def\BB{2}
\def\CC{3}

\def\II{9}

\font\titfont=cmr10 at 12 pt

\def\f{\varphi}
\def\z{\zeta}

 \def\II{1}
 \def\AA{2}
 \def\BB{3}
\def\CC{4}
   
\  \vskip .2in
\centerline{\titfont  THE  PROJECTIVE HULL OF CERTAIN CURVES IN C$^2$.  }
 
\vskip .2in
\centerline{\titfont Reese Harvey, Blaine Lawson$^*$ and John Wermer}

\vglue .9cm
\smallbreak\footnote{}{ $ {} \sp{ *}{\rm Partially}$  supported by
the N.S.F. }

\vskip .3in
\centerline{\bf Abstract} \medskip
  \font\abstractfont=cmr10 at 10 pt

{{
\parindent= .7in\narrower\abstractfont \noindent
The projective hull $\PH X$ of a compact set $X\subset\bbp^n$ is an 
analogue of the classical polynomial hull of a  set in $\bbc^n$. In the 
special case that $X\subset \bbc^n\subset \bbp^n$, the affine part
$\PH X \cap \bbc^n$ can be defined as the
set of points $x\in\bbc^n$ for which there exists a constant $M_x$
so that

}}
$$
|p(x)|\ \leq \ M_x^d\sup_X|p|
$$
{{\parindent= .7in\narrower\abstractfont \noindent for all polynomials $p$ of degree $\leq d$, and any $d\geq 1$. 
Let $\PH X (M)$ be the set of points $x$ where $M_x$ can be chosen $\leq M$.
Using an argument of E. Bishop, we show that if $\g\subset \bbc^2$ is a compact real analytic curve
(not necessarily connected), then for any linear projection $\pi:\bbc^2\to\bbc^1$, 
the set {\sl $\PH \g (M)\cap \pi^{-1}(z)$ is finite for almost all $z\in\bbc$}. It is then shown that
for any compact stable  real-analytic curve $\g\subset\bbp^n$, the set {\sl $\PH \g-\g$ is a 
1-dimensional complex analytic subvariety of $\bbp^n-\g$}.

}}

\vfill\eject

\noindent
{\bf \S\II. Introduction.} 
The classical {\bf polynomial hull} of a compact subset $X\subset \bbc^n$ is the set of points
$x\in \bbc^n$ such that   \def\s{P}
$$
|p(x)|\ \leq \ \sup_X |p|\qquad {\rm for\  all\  polynomials \ }\ p.
\eqno{(\II.1)}
$$
In [4] the first two authors introduced an analogue for compact subsets of projective space.
Given $X\subset \bbp^n$, the {\bf projective hull of $X$} is the set $\PH X$ of points
$x\in \bbp^n$ for which there exists a constant $C=C_x$ such that 
$$
\|\s(x)\|\ \leq\ C_x^d \sup_X \|\s\|   \qquad {\rm for\  all\    sections \ }\ \s\in H^0(\bbp^n,\co(d)) 
\eqno{(\II.2)}
$$
and all $d\geq 1$. Here $\co(d)$ is the $d$th power of the hyperplane bundle with its 
 standard metric. Recall that $H^0(\bbp^n,\co(d))$ is given naturally as the set of homogeneous 
 polynomials of degree $d$ in homogeneous coordinates.  If $X$ is contained in an affine chart
 $X\subset \bbc^n\subset \bbp^n$ and $x\in\bbc^n$, then condition (\II.2) is equivalent to
   \def\s{\sigma}
  $$
|p(x)|\ \leq \   M_x^d\sup_X |p|\qquad {\rm for\  all\  polynomials \ }\ p  {\rm \ of \ degree\ } d
\eqno{(\II.3)}
$$
  and all $d\geq 1$ where $M_x =   \rho\sqrt{1+\|x\|^2} C_x$ and $\rho$ depends only on
  $X$.  Therefore the set $\PH X \cap \bbc^n$ consists exactly of those points $x\in \bbc^n$
  for which there exists an $M_x$ satisfying condition  (\II.3).
  
  This paper is concerned with the case where $X=\g$ is a real analytic curve. In [4] 
  evidence was given for the following conjecture.
  
  \Conj{\II.1} {\sl Let $\g \subset \bbp^n$ be a finite union of simple closed real analytic curves.
  Then $\PH \g-\g$  is a 1-dimensional complex analytic suvariety  of  $\bbp^n-\g$.}
  
  \medskip
  This conjecture
   has many interesting geometric consequences (See [5], [6], and [7] ).
   
  The assumption of real analyticity is important. The conjecture does not hold 
   for all smooth curves.  In particular, it does not hold for curves which are not pluripolar.

   One point of this paper is to prove Conjecture \II.1 under the hypothesis that the 
   function $C_x$ is bounded on $\PH \g$.
We  begin by adapting arguments of E. Bishop in [2] to prove the following finiteness theorem.
  
  \Theorem{\II.1} {\sl Let $\g \subset \bbc^2$ be a finite union of simple closed real analytic curves.
  Set
  $$
 { \PH\g}_M\ \equiv \ \{x \in \PH\g\cap\bbc^2 : M_x\leq M\}
  $$
where $M_x$ is the function appearing in  condition (\II.3). Let $\pi:\bbc^2\to \bbc$ be a linear projection.
Then 
$$
{ \PH\g}_M\cap \pi^{-1}(z) \ \ {\sl is\ finite\ for \ almost\ all\ } z\in\bbc.
$$
Consequently, ${ \PH\g}\cap \pi^{-1}(z)$ is countable for almost all $z\in\bbc$.
}
  \medskip
  
  In section 3 this theorem is combined with results from [4] and the theorems concerning
  maximum modulus algebras to prove the following. A set $X\subset \bbp^n$ is called {\sl stable}
  if the function $C_x$ in (\II.2) is bounded on $\PH X$.  Note that if $X$ is stable and $X\subset \bbc^n\subset \bbp^n$, then the function $M_x$ is bounded on $\bbc^n$ by $\rho\sqrt{1+\|x\|^2}$.

  \Theorem{\II.2} {\sl Let $\g \subset \bbp^n$ be a finite union of simple closed real analytic curves.
  Assume $\g$ is stable.  Then $\PH \g-\g$  is a 1-dimensional complex analytic subvariety  of  $\bbp^n-\g$.}

  \vfill\eject
  

\noindent
{\bf \S\AA. The Finiteness Theorem.} 
Let $X$ be a compact set in $\bbc^n$ and denote by $\cp_d$
 the space of polnomials of degree $\leq d$ on $\bbc^n$.

 \Def{\AA.1} Denote by $\wh X \cap \bbc^n$ the set of all $x\in \bbc^n$ such that there exists a constant
 $M_x$ with 
 $$
 |P(x)|\ \leq \ M_x^d\cdot \sup_X|P|
 \eqno{(\AA.1)}
$$ 
 for every $P\in\cp_d$ and  $d\geq1$. The set  $\wh X\cap \bbc^n$ is called {\bf the projective hull of $X$ in } $\bbc^n$.

 \medskip
 As noted above, the projective hull,  defined in [4],  is a subset of projective space $\bbp^n$,
 and the set $\wh X\cap \bbc^n$    is exactly that part of the projective hull which lies in the  affine chart
 $\bbc^n\subset\bbp^n$.  
Closely related to Definition \AA.1  is the following.

 \Def{\AA.2} Fix a number $M\geq 1$ and a point $z\in\bbc^{n-1}$.  Then we set
 $$
 {\wh X}_M(z)\ =\  \{w\in \bbc : |P(z,w)| \leq M^d\cdot \sup_X|P| \ \forall P\in\cp_d \ {\rm and }\ \forall d\geq 1\}
 $$
and let $\wh X(z) =\bigcup _{M\geq1}  {\wh X}_M(z)= \{ w\in\bbc : (z,w)\in  \wh X\}$.
 
 \medskip
 We consider a special case of these definitions.   We fix $n=2$ and consider a simple closed 
 real-analytic curve $X$ in $\bbc^2$.  Let $\D$ denote the unit disk in $\bbc$.
 
 \Theorem{\AA.1}  {\sl  Fix $M\geq1$.  For almost all $z\in\D$, $ {\wh X}_M(z)$ is a finite set.}

\Cor{\AA.1} {\sl For almost all $z\in \bbc$ the set $ {\wh X}(z)$ is countable.}\medskip
 
 We shall prove Theorem \AA.1 by adapting an argument,  for the case of polynomially convex hulls,
 by Errett Bishop in [2].  We shall follow the exposition of Bishop's argument in [10], Chapter 12.
 
 \Def{\AA.3} The polynomial
 $
 Q(z,w)\ =\ \sum_{n,m} c_{nm} z^n w^m
 $
 is called a {\bf unit polynomial} if $\max_{n,m} |c_{nm}| =1$.
  
  \Def{\AA.4}  The polynomial
  $
 Q(z,w)\ =\ \sum_{n,m} c_{nm} z^n w^m
 $
 is said to have  {\bf bidegree} $(d,e)$, for non-negative integers $d$ and $e$,  if $c_{nm}=0$ unless $n\leq d$ and $m\leq e$.
  
  \medskip Note that $\deg\, Q \leq d+e\leq 2\deg\, Q$.
  
  \Def{\AA.5}  Fix $M\geq 1$.  For each $z\in\bbc$ set
  $$
  S_M(z)\ =\ \{w\in \bbc : |Q(z,w)|\leq(M^{d+e})\sup_X|Q| \ \ \forall Q \in \bbc[z,w] \ 
  {\rm of \  bidegree\ }  (d,e) \ {\rm for\ } d,e\geq 1\}.
  $$
  
  We now fix a number $M\geq 1$ and keep it fixed throughout what follows.

  \Theorem{\AA.2}  {\sl For almost all $z\in \D$, $S_M(z)$ is a finite set.}
  \medskip
  
  Theorem \AA.1 is an immediate consequence of Theorem \AA.2.  To see this, fix $z\in \D$ and choose $w\in  {\wh X}_M(z)$. Choose next a polynomial $Q$ of bidegree $(d,e)$ and let $\d=\deg\, Q$. Then
  $$
  |Q(z,w)|\ \leq\ M^\d\| Q\|_X  \ \leq\ M^{d+e}\|Q\|_X
  $$
  and so $w\in S_M(z)$. Since this holds for all such $w$, ${\wh X}_M(z)\subseteq S_M(z)$.
  By Theorem \AA.2 $S_M(z)$ is a finite set for a. a. $z\in \D$.  so ${\wh X}_M(z)$ is a finite set
  for a. a. $z\in \D$. Thus Theorem \AA.1 holds.
  
  We now go to the proof of Theorem \AA.2.
  
  \Lemma{\AA.1}  {\sl Let $\O$ be a plane domain, let $K$ be a compact set in $\O$, and fix $z_0\in\O$.
  Then there exists a constant $r$, $0<r<1$, so that if $f$ is holomorphic on $\O$ and $|f|<1$ on $\O$
  and if $f$ vanishes to order $\l$ at $z_0$, then }
  $$
  |f|\ \leq \ r^\l\qquad {\sl on}\ K.
  $$

  \pf  We construct a bounded and smoothly bounded subdomain $\O_0$ of $\O$ with
   $\overline{\O_0} \subset \O$, $z_0\in \O_0$ and $K\subset \O_0$.
  Denote by $G(z_0,z)$ the Green's function of $\O_0$ with pole at $z_0$.
  
  Then $e^{-(G+iH)}$ is a multiple-valued holomorphic function on $\O_0$ with a 
  single-valued modulus $e^{-G}$, and this modulus is  $=1$ on $\bo_0$.   ($H$ is the harmonic conjugate of $G$.)  Consequently, $$f \over  e^{-\l(G+iH)}$$ is multiple-valued and holomorphic
  on $\O_0$, and  its modulus is single-valued and $<1$ on $\bo_0$.
  By the maximum principle for holomorphic functions, for each $z\in K$, we have
  $$
  \left|  { f \over e^{-\l(G+iH)} }\right|\ < \ 1
  $$
   at $z$ and so 
  $$
  |f(z)
   \ \leq \ \left[ e^{-G(z_0,z)}\right]^\l.
  $$
  Putting $r=\sup_K e^{-G}$, we get our desired inequality.\qed
  
  \Lemma{\AA.2}   {\sl Let $\O$ be a bounded plane domain and $K$ a compact subset of $\O$.
  Let $\cl$ be an algebra of holomorphic functions on $\O$.  Put $\|\f\|=\sup_K|\f|$
for all $\f\in\cl$.
  
  Fix $f, g \in \cl$.  Then there exist $r$, $0<r<1$ and $C>0$ such that for each pair of positive 
  integers $(d,e)$ we can find a unit polynomial $F_{d,e}$ of bidegree $(d,e)$ such that }
  $$
  \left\| F_{d,e}(f,g)\right\|\ \leq \ C^{d+e}r^{de}.
  \eqno{(\AA.2)}
  $$
  
  \pf
  Choose a subdomain $\O_1$ of $\O$ with $K\subset \O_1 \subset  \overline {\O_1}\subset
  \O$.  Choose $C_0>1$ with $|f|<C_0$, $|g|<C_0$ on $\overline {\O_1}$. 
  Consider an arbitrary polynomial
  $$
  F(z,w)\ =\ \sum_{n=0}^d\sum_{m=0}^e c_{nm} z^n w^m
  $$
  and let $h$ be the function $F(f,g)$ in $\cl$.  
  Fix a positive integer $\l$. The requirement that $h$ should vanish at $z_0$ to order $\l$
  imposes $\l$ linear homogeneous conditions on the $c_{nm}$, 
  and hence has a non-trivial solution  if $\l <(d+1)(e+1)$. We may assume that the corresponding polynomial $F$ is a unit polynomial.  Since
  $$
  { d^\nu h \over dz^\nu }(z_0) \ =\ 0, \ \ \ \nu=0,1,...,\l-1,
  $$
  Lemma \AA.1 gives us some $r$, $0<r<1$, such that
  $$
  |h|\ \leq\ \left( \sup_{ \overline {\O_1}}  |h|  \right) \cdot r^\l \qquad{\rm on\ \ } K.
  $$
  Since $F$ is a unit polynomial,
  $$
  |h| \ \leq \ \sum_{n=0}^d\sum_{m=0}^e  |c_{nm}| |f|^n |g|^m\ \leq\ (d+1)(e+1) C_0^{d+e} 
  \qquad{\rm on\ \ }  \overline {\O_1}.
  $$
Hence for large $C$,
$$
\|h\|\ \leq\  (d+1)(e+1) C_0^{d+e} \ \leq\ C^{d+e} r^\l.
$$
  We choose $\l=de$.  Since $de< (d+1)(e+1)$, we get
  $$
  \|F(f,g)\| \ =\ \|h\|\ <\ C^{d+e} r^{de}
  $$
  as desired.\qed
  
  \Note{}  We shall  apply this result to the case when $K$ is the unit circle, 
  $\O$ is an annulus containing $K$, and $\cl$  is the algebra of functions
  holomorphic on $\O$.
  \medskip
  
  The curve $X$ in our Theorem \AA.2 is real analytic by hypothesis, and hence can
  be represented parametrically:
  $$
  z\ =\ f(\z), \ \ w\ =\ g(\z) \qquad\z\in\O
  $$
  where $f,g$  are functions in $\cl$.
  
  \Lemma{\AA.3}  {\sl Let $r,C$ and $F_{d,e}$ be as in Lemma \AA.2.  Fix $r_0$, $r<r_0<1$.
  Then there exists $d_0$ such that}
  $$
  (MC)^{d+e}\cdot r^{de}\ \leq\ r_0^{de} \qquad {\sl for\ \ }d,e >d_0.
  \eqno{(\AA.3)}
  $$
  \pf
  We write $\sim$ for ``is equivalent to''.
  $$ \eqalign{
  (\AA.3)\ \ \ &\sim\ \ (MC)^{d+e}\leq \left({r_0\over r}\right)^{de}    \cr
  &\sim\ \ (d+e)\log(MC) \leq  de \log\left({r_0\over r}\right)    \cr
  &\sim\ \ \left({1\over e}+{1\over d}\right)\log(MC) \leq    \log\left({r_0\over r}\right).    \cr
  }
  $$
  The last inequality is true for $d,e>d_0$ for some suitable $d_0$.  We are done.\qed
  
  \medskip
  
  With $M,r,r_0$ fixed, we choose $d_0$ as in (\AA.3).  Henceforth, we tacitly assume $d,e>d_0$.
  
  \Def{\AA.6}  Fix $d,e$ and put $F=F_{d,e}$  as above.
  Then 
  $$
  F(z,w)\ =\ \sum_{j=0}^e G_j(z) w^j
  $$
  where for some $j=j_0$, $G_{j_0}$ is a unit polynomial of degree $ \leq  d$.
  We define
  $$
  T(d,e)\ =\   \left\{z\in \D : \left| G_{j_0} (z)\right| \ \leq\ r_0^{de\over 2}\right\}.
  $$
  
  \Lemma{\AA.4}  {\sl Let $F$ be a unit polynomial in $z$, of degree $k$, and let 
  $\a$ be a positive number.  Put $\L = \{z\in\D : |f(z)|\leq \a^k\}$. Then
  $$
  m(\L)\ \leq\ 48 \a,
  $$
  where $m$ is two-dimensional measure.}
  
  \pf  This is Lemma 12.3 in [10], and a proof of it is given there.
  
  \Lemma{\AA.5}  {\sl Fix $d,e$. Fix a point $z_1\in \D- T(d,e)$.  Then there exists a unit polynomial
  $B$ in one variable, of degree $\leq e$, such that for every $w_0\in S_M(z_1)$, we have }
  $$
  |B(w_0)|\ \leq \ r_0^{de\over 2}.
    $$
\pf  
Define the polynomial $A$ in one variable by $A(w) = F(z_1,w)$, where $F=F_{d,e}$.  As in Definition  
  \AA.6 then
  $$
  A(w)\ =\  \sum_{j=0}^e G_j(z_1) w^j
  $$
  and $G_{j_0}$ is a unit polynomial in $z$. 
  
   Since  $z_1 \notin T(d,e)$, we have
   $$
   \left|  G_{j_0}(z_1)   \right|\ >\ r_0^{de\over 2}.
   \eqno{(\AA.4)}
   $$  
  Fix $w_0 \in S_M(z_1)$.  Then
  $$\eqalign{
  |F(z_1,w_0)|\ &\leq\ M^{d+e}\cdot \| F\|_X  \cr
  &\leq\ M^{d+e} C^{d+e}\cdot r^{de} \qquad \qquad {\rm by\ \ } (\AA.2) \cr
 &\leq\  r_0^{de} \qquad \qquad\qquad \qquad \ \ \ \ \ {\rm by\ \ } (\AA.3).\cr
  }
  $$
  We shall divide $A$ by its largest coefficient $K$. Note that
  $$|K|\geq | G_{j_0}(z_1) | > r_0^{de\over 2}$$
  by (\AA.4).  Put $B(w) = A(w)/K$.  Then deg$\,B \leq e$
  and
  $$
  |B(w_0)|\ =\ {|A(w_0)|\over |K|} \ = \  {|F(z_1,w_0)|\over |K|} \ \leq \ { r_0^{de}  \over   r_0^{de\over 2} }
  \ =\   r_0^{de\over 2}.
  $$
  We are done.\qed
  
  \Lemma{\AA.6}  {\sl For each $d$,}
  $$
  m(T(d,e)) \ \leq \ 48\, r_0^{e\over 2}
  $$
  \pf
  Fix $e$ and fix $d$.  With $G_{j_0}$ as above, write $G=G_{j_0}$.   Then deg$\,G\leq d$.
  By definition of $T(d,e)$, if $z\in T(d,e)$, then
  $$
  |G(z)\ \leq \ r_0^{de\over 2}\ =\ \left(r_0^{e\over 2}\right)^d \ \leq\ \left(r_0^{e\over 2}\right)^{\deg\, G},
  $$
 and  so
  $$
  T(d,e)\ \subseteq \ \left\{z\in\D : |G(z)| \ \leq \ \left(r_0^{e\over 2}\right)^{\deg\, G}\right\}.
  $$
  Therefore,
  $$
 m\left[ T(d,e)\right] \ \leq\ m\left\{  z\in\D : |G(z)| \ \leq \  \a^k  \right\}
  $$
  where $\a = r_0^{e\over 2}$ and $k=\deg\,G$.  By Lemma \AA.4,
  $m\left\{  z\in\D : |G(z)| \ \leq \  \a^k  \right\} \leq 48\,\a$, and so 
  $ m\left[ T(d,e)\right] \leq 48\,  r_0^{e\over 2}$, as was to be shown.\qed
  
  \Def{\AA.7}  Fix $e$ and and set
  $$
  H_e 
   =\ \{z : z\in \D-T(d,e)\ \ {\rm for\ infinitely \ many\  } d\}.
  $$
  
  \Lemma{\AA.7}  {\sl If $z^* \in H_e$, then $S_M(z^*)$ has at most $e$ elements.}
  
  \pf
  Fix $z^* \in H_e$. Then there exists a sequence $\{d_j\}$ such that $z^*\in \D-T(d_j,e)$ for each $j$.
  By Lemma \AA.5, for each $j$ there is a unit polynomial $B_j$ with $\deg\, B_j\leq e$ such that
  $$
  \left|B_j(w_0)\right|\ \leq\ r_0^{d_j e\over 2} \qquad\ \ {\rm for\ each\ \ }w_0\in S_M(z^*).
  \eqno{(\AA.5)}
  $$
  
  Since$\deg\,B_j\leq e$ for all $j$, and each $B_j$ is a unit polynomial, there exists a subsequence of the sequence $\{B_j\}$ converging uniformly to a unit polynomial $B^*$ on compact sets in the $w$-plane.  Because of (\AA.5), $B^*(w_0)=0$ for each $w_0\in S_M(z^*)$. Also,  $\deg\, B^*\leq e$.
  Hence the cardinality of $S_M(z^*)$ is $\leq e$. We are done.\qed
  \medskip
  
  \noindent
  {\bf Proof of Theorem \AA.2}.  Our task is to show that 
  $m\{z\in\D: S_M(z) $ is infinite $\}=0$.  
  Fix $e$.  Fix $z\in\D-H_e$.  Since $z\notin H_e$, we have $z\in \D-T(d,e)$ for only finitely many $d$, so
  $z\in T(d,e)$ for all $d$ from some $d=k$ on.  Therefore,
  $$
  z\ \in \ \bigcap_{d=k}^\infty T(d,e)
  $$
  and so
  $$
  \D -H_e \ \subseteq\ \bigcup_{d=d_0}^\infty\left[ \bigcap_{d=k}^\infty T(d,e)\right].
  \eqno{(\AA.6)}
  $$
  
  By Lemma \AA.6, $m(T(d,e) \leq 48\, r_0^{e\over 2}$ for each $d$.  Therefore,
  $$
  m\left(\cap_{d=k}^\infty T(d,e)\right) \ \leq 48\, r_0^{e\over 2}
  $$
  for each $k$. So the right hand side of (\AA.6) is the union of an increasing family
  of sets each of which has $m$-measure $\leq 48\, r_0^{e\over 2}$.  Thus (\AA.6) gives
  $$
  m\left(\D-H_e\right)\ \leq\ 48\, r_0^{e\over 2}.
  \eqno{(\AA.7)}
  $$
  Also, by Lemma \AA.7, we have
   $$
  {\rm If }\ \ z^*\in H_e, \ {\rm then\ }\#\left\{ \,S_M(z^*)\right\} \ \leq \ e.
  \eqno{(\AA.8)}
  $$
  Fix $z\in\D$ such that the set $S_M(z)$ is infinite.
  Then $z\notin H_e$ for each $e$, that is, $z\in\D-H_e$ for all $e$.
  Hence,  $\{z\in\D: S_M(z) $ is infinite $\} \subset \D-H_e$.  Therefore
   $$ 
   m\{z\in\D: S_M(z) \ \ {\rm  is \ infinite\ } \} \ \leq \ m(\D-H_e)\leq 48\, r_0^{e\over 2}
   $$
   by (\AA.7).  We now let $e\to\infty$ and conclude that 
   $m\{z\in\D: S_M(z) $ is infinite $\}=0$.  Theorem \AA.2 is proved. \qed
  
  \medskip
  \noindent
  {\bf Proof of Corollary \AA.1}. Fix $r>0$ and  apply Theorem \AA.1 to the curve $\rho_r(X)$
 where $\rho_r:\bbc^2\to\bbc^2$  is given by $\rho_r(z)=rz$.
  Since   $\rho_r(\PH X\cap\bbc^2) =\PH{(\rho_r X)}\cap\bbc^2$, we conclude that Theorem \AA.1 holds  with $\D$ replaced by
 $\smfrac 1 r \D$.\qed
  
  \Theorem{\AA.3} {\sl Theorem \AA.1 remains valid without the assumption that $X$ is connected,
  that is, it is valid when $X$ is a finite union of real analytic simple closed curves in $\bbc^2$.}
  
  \pf
  Write $X = \g_1\cup\g_2\cup\cdots\cup\g_N$ where each $\g_k\subset \bbc^2$ is a simple closed
  real analytic curve.  Choose a neighborhood $\O$ of the unit circle $K$ in $\bbc$ and   
  complex analytic maps $(f_k,g_k):\O_k\to\bbc^2$, $k=1,...,N$ whose  restriction to $K$ is
  a parameterization of $\g_k$.    We now apply the following.

  \Lemma{\AA.8}   {\sl Let $\O$ be a plane domain and $K$ a compact subset of $\O$.
  Let $\cl$ be an algebra of holomorphic functions on $\O$.  Put $\|\f\|=\sup_K|\f|$
for all $\f\in\cl$.
  
  Fix $f_k, g_k \in \cl$ for $k=1,...,N$.  Then there exist $r$, $0<r<1$ and $C>0$ such that for each pair of positive 
  integers $(d,e)$ with $d+e>N$, we can find a unit polynomial $F_{d,e}$ of bidegree $(d,e)$ such that }
  $$
  \left\| F_{d,e}(f_k,g_k)\right\|\ \leq \ C^{d+e}r^{de\over N}\qquad{\rm for\ \ } k=1,...,N.
  \eqno{(\AA.9)}
  $$
\pf We fix a point $z_0\in \O$ and choose $F_{d,e}$ so that $F_{d,e}(f_k,g_k)$ vanishes to order
$de/N$ at $z_0$ for all $k$.  This is possible if $d+e>N$.  We then proceed as in the proof of 
Lemma \AA.2.\qed
\medskip

One can now carry out the arguments given above for the case of one component.
The only difference is that in the estimates, $r_0^e$ will be replaced by $r_0^{e\over N}$. \qed


 \vfill\eject
 
 \noindent
 {\bf \BB.\ The Analyticity Theorem.}  Let  $\co(1)\to \bbp^n$ denote 
 the holomorphic line bundle of Chern class 1 over complex projective $n$-space, endowed 
 with its standard U(n+1)-invariant metric $\|\cdot\|$. Following [4], we define the  {\sl projective hull} of  
 a compact subset $X\subset \bbp^n$  to be the set $\PH X$ of points  $x\in \bbp^n$ for which 
 there exists a constant $C=C_x$ such that
 $$
 \|\s(x)\|\ \leq \    C_x^d \sup_X \|\s\|.
\eqno{(\BB.1)} $$
 for all holomorphic sections $\s\in H^0(\bbp^n,\co(d))$ and all $d\geq 1$.

 \Note{\BB.1} Recall that
the  holomorphic sections $H^0(\bbp^n,\co(d))$ correspond naturally to the homogeneous
polynomials of degree d in homogeneous coordinates $[Z_0,...,Z_n]$ for $\bbp^n$. From this 
one can see (cf. [4, \S 6])  that if $X$ is contained in an affine chart $\bbc^n\subset\bbp^n$, then
$\PH X\cap\bbc^n$ is exactly the ``projective hull of $X$ in $\bbc^n$'' introduced in \S \AA. 
Moreover, the function $M_\z$ appearing in  (\AA.1) can be taken to be $M_\z = \rho\sqrt{1+\|\z\|^2}C_\z$ for  $\z\in \PH X \cap\bbc^n$, where $\rho$ is a constant depending only on $X$.
   \medskip
   
   For each $x\in\PH X$ there is a {\sl best constant} 
   $C(x) \equiv \min\{C_x : $   (\BB.1) holds   $\forall\s \}$. The set $X$ is called {\bf stable} if the best constant function $C$ is bounded  on $\PH X$.
   We know from [4, Prop. 10.2] that if $X$ is stable, then $\PH X$ is compact.

   The point of this section is to prove the following projective version of the main theorem
   in [9].
   
   \Theorem{\BB.1}  {\sl  Let $\g\subset \bbp^n$ be a finite union of real analytic closed
   curves and assume $\g$ is stable.  Then $\PH \g-\g$ is a one-dimensional complex 
   analytic subvariety of $\bbp^n-\g$.}
      
  \Note{\BB.2}  
  When  this conclusion holds,  one can show that, in fact, $\PH \g$ is the image of 
   a compact  riemann surface with analytic boundary 
   under a holomorphic map to $\bbp^n$.  We will prove this  in \S \CC.
  
  \pf  Assume to begin that $n=2$.
  Since $\g$ is real analytic, it is pluripolar, i.e., locally contained in the 
  $\{-\infty\}$-set of a plurisubharmonic function  (which is $\not \equiv -\infty$). 
  Therefore, by [4, Cor. 4.4]  we know that $\PH \g$ is also pluripolar.  In particular,
  it is nowhere dense.  As noted above,  $\PH \g$ is closed by stability.  Hence, we may choose
  a point $x\in \bbp^2$ and a ball $B$ centered at $x$ such that 
  $$
\PH\g  \ \subset\ \bbp^2- \overline B.
   $$
  Let 
  $$
\bbp^2-\{x\} \     \Arr{\pi}    \ \bbp^1
  \eqno{(\BB.2)}
   $$
  be linear projection  with center $x$.
  This projection  (\BB.2) is naturally  a holomorphic  line bundle   $\cong \co(1)$,  and 
  $$
  \bbp^2-\overline B      \Arr{\pi}    \ \bbp^1
  \eqno{(\BB.3)}
    $$
  can be identified, after scalar multiplication by some constant $r>0$,
  with its  open unit disk bundle.
  
  \def\U{V}
  
 Cover $\bbp^1$ with two affine charts: $\U_0 = \bbp^1-\{0\}$ and 
  $\U_{\infty} = \bbp^1-\{\infty\}$, and  assume that $\g\cap\pi^{-1}(0) =
  \g\cap\pi^{-1}(\infty) =\emptyset$. By symmetry we may  restrict attention
  to  
  $
  \pi^{-1}(\U_\infty)
  $.
   This chart has an   identification 
    $$
  \pi^{-1}(\U_\infty) \cong \bbc^2 =\{(z,w): z,w\in\bbc\}
  $$
  with the property that $\U_\infty$ maps linearly to the $z$-axis and 
  $\pi$ can be written as $\pi(z,w)=z$.
  The subset $ \bbp^2-\overline B$, intersected with this chart, is represented by
  $$
   (\bbp^2-\overline B)\cap \bbc^2\ =\ \{(z,w) : |w|^2 \leq  |z|^2+1\}.  
 \eqno{(\BB.4)}
    $$
  Set 
  $$
  \O \equiv \bbc-\pi(\g)
  \qquad{\rm and }\qquad
  U\equiv\pi^{-1}(\O) = \bbc^2-\pi^{-1}(\pi(\g)).
  $$
  \Prop{\BB.1}  {\sl Let $\g\subset \bbc^2$ be a stable real analytic curve
  with the property that 
  $$
   \PH \g  \cap \bbc^2 \ \subset \ \{(z,w) : |w|^2 \leq  |z|^2+1\}.
  \eqno{(\BB.5)}
    $$
  Then $ \PH \g  \cap U$ is a 1-dimensional complex analytic subvariety   of   $U$.}

 \pf
  Note to begin that since $\PH \g$ is compact, condition (\BB.5) implies that
  $$
  \pi:\PH \g \cap U\ \to \O \qquad {\rm is\ a\  proper\ map}.
  \eqno{(\BB.6)}
  $$
  Consider now the algebra $A$ of functions on $\PH \g\cap U$ given by restriction of the 
  holomorphic functions on $U$, i.e., 
  $$
  A\ \equiv\     \left\{f\bigr|_{\PH \g \cap U} : f\in \co{}(U)   \right\}.
  $$
  We now claim that $(A, \PH \g \cap U, \O, \pi)$ is a {\sl maximum modulus algebra}, as defined
  in [1, pg.64].  Given (\BB.6) this means that we need only prove the following.
  
  \Lemma {\BB.1}  {\sl For each $z_0\in \O$ and each closed disk $D\subset\O$ 
  centered at $z_0$, the equality
  $$
  |f(z_0, w_0)|\ \leq \ \sup_{\PH \g \cap \pi^{-1}(\partial D)} |f|
\eqno{(\BB.7)}
  $$
  holds for all $f\in A$.}
  
  \pf
  By hypothesis (\BB.5) there exists an $R>0$ such that 
  $$
 \PH \g \cap \pi^{-1}(D)  \ \subset \ D\times \D_{R/2}
  $$
  where $\D_r \equiv \{w : |w| \leq r\}$.  In particular,  we have that
 $$
  \PH \g \cap \partial (D\times \D_R) \ =\  \PH \g \cap (\partial  D\times \D_R)
    \ =\ \PH \g \cap \pi^{-1}(\partial D).
 \eqno{(\BB.8)}
  $$
Now   Theorem 12.8 in [4]  states that 
  $$
   \PH \g \cap \pi^{-1}(D) \ =\ 
  \PH \g \cap (D\times \D_R)\ \subset \ {\rm Polynomial\ Hull\ of\ } \PH \g \cap \partial (D\times \D_R).
  $$
 Applying  (\BB.8)  gives
    $$
 \PH \g \cap \pi^{-1}(D) \ \subset \ {\rm Polynomial\ Hull\ of\ } \PH \g \cap \pi^{-1}(\partial  D),
   $$
and Lemma \BB.1 follows immediately.
\qed
\medskip

We have now shown that $(A, \PH \g \cap U, \O, \pi)$ is a   maximum modulus algebra.
Furthermore, since $\PH\g$ is stable, we know from Theorem \AA.1 that
  there exists an $N>0$ such that
  $$
  \O(N) \ \equiv \   \left \{z\in\O : \#\left(\pi^{-1}(z) \cap \PH\g\right)\ \leq \ N\right\}
  $$
  has positive measure. (Since $\O-\bigcup_N\O(N)$ has measure zero.)
  It now follows from Theorem 11.8 in [1] that:
  \smallskip
  
   \item{(i)}  $\O=\O(N)$, and  \smallskip
  
   \item{(ii)}  There exists a discrete subset $\L\subset \O$ such that $\PH \g \cap\pi^{-1}(\O-\L)$
 has the structure of a Riemann surface on which every function in $A$ is analytic.
 \smallskip

 \noindent
 Since $A$ is the restriction of holomorphic functions on $U$ to $\PH\g$, condition (ii) implies that
 $\PH\g \cap \pi^{-1}(\O-\L)$ is a 1-dimensional complex analytic subvariety of $\pi^{-1}(\O-\L)
 = U-\pi^{-1}(\L)$.

 It now follows that  $\PH\g \cap U$ 
 is a 1-dimensional complex analytic subvariety of $U$. To see this, fix $z_0\in\L$
 and choose a small closed disk $D\subset \O$ centered at $z_0$ with $D\cap \L=\emptyset$.
 The arguments above show that $\PH\g\cap \pi^{-1}(D)$ is contained in the polynomial 
 hull of the real analytic curve $\PH\g\cap \pi^{-1}(\partial D)$.  Applying 
 standard  results [1, \S 12] 
 proves Proposition \BB.1  \qed
 \medskip
 
 Proposition \BB.1 together with the discussion preceding it,  give the following.
 
  \Cor{\BB.1}  {\sl The set $\PH\g-\pi^{-1}(\pi\g)$ is a 
complex analytic subvariety   of dimension one in  $\bbp^2-\pi^{-1}(\pi\g)$.}
\medskip

 Observe  that for every point   $y\in \bbp^2-\PH \g$ there is a point $x\in \bbp^2-\PH \g$
 such that $\pi(y)\notin\pi(\g)$ where $\pi$ is the projection (\BB.2) with center $x$.
 Consequently, Corollary \BB.1 proves   Theorem \BB.1 for the case $n=2$.
 
 Suppose now that $n=3$ and choose $x\in \bbp^3-\PH\g$. 
  The set of such $x$ is open and dense since 
 $\PH\g$ is a compact pluripolar set of Hausdorff dimension 2 (cf.  [4, Cor. 4.4 and Thm. 12.5]).  
 Let $\Pi:\bbp^3-\{x\}\to\bbp^2$
 be the projection with center $x$. 
One sees easily that $$\Pi(\PH \g)\ \subseteq\  \PH{\Pi\g},$$
and by the above $\PH{\Pi\g}-{\Pi\g}$ is a complex analytic curve in $\bbp^2-{\Pi\g}$.
  Standard arguments  now show that $\PH \g-\g$ is a complex analytic curve
  in $\bbp^3-\g$. Proceeding by induction on $n$ completes the proof of Theorem \BB.1.\qed

    \vfill\eject
 \noindent
 {\bf \CC. Boundary Regularity.}   The conclusion of Theorem \BB.1 implies a strong regularity
 at the boundary. For  future reference
 we include a discussion of this regularity.
 
 \Theorem{\CC.1}  {\sl Let $\g\ss \pn$ be a finite union of real analytic closed curves,
 and suppose $V$ is a 1-dimensional complex analytic subvariety of the complement
 $\pn-\g$.  Then 
 $$
 V\ =\ \bigcup_{j=1}^m V_j \,\cup\, \bigcup_{k=m+1}^\ell V_k'  \qquad {\sl where}
 $$
 (1)\ \ Each $V_j$ is an irreducible 1-dimensional complex analytic subvariety
 of finite area in $\pn-\g$ whose closure  $\overline{V}_j$ is an immersed variety in $\pn$ 
  with non-empty
 boundary $\partial \overline{V}_j = \g_j$ consisting of a union of components of $\g$.
 In particular, there exists a connected Riemann surface $S_j$, a compact subdomain
 $\overline{W}_j\ss S_j$ with real analytic boundary, and a   generically injective holomorphic map
 $$
 \rho_j: S_j\ \arr\ \pn \qquad{\sl with} \ \ \rho_j(\overline{W}_j)\ =\ \overline{V}_j
 $$
 which is an embedding on a neighborhood of $\partial \overline{W}_j$ and has 
 $\rho_j(\partial \overline{W}_j)=\g_j$.
 \medskip
 \noindent
 (2)\ \ Each $V_k'$ is an irreducible algebraic curve in $\pn$ with $\g_k\ss  {\rm Reg}(V_k')$ where 
 $\g_k$ is a (possibly empty) finite union of components of $\g$.
 \medskip
 \noindent
 (3) \ \  The curve $\g$ is a disjoint union $\g=\g_0\cup \g_1\cup\cdots\cup \g_\ell$
 where $\g_0$ is also a finite union of connected components of $\g$.
 }

 \Note{\CC.2}  When $\g$ is stable and $V=\wh\g$, each $\g_k$ is non-empty for $m<k\leq \ell$.

  \medskip
  Theorem \CC.1 can be put into a more succinct form.
  
  \Theorem {\CC.1$'$} {\sl Let $\g$ and $V$ be as above. Then there exists a Riemann surface
  $S$ (not necessarily connected), a compact subdomain $\overline W\ss S$ with real analytic boundary, and a holomorphic
  map $\rho:S\to \pn$ which is generically injective and satisfies
  \medskip\noindent
  (1) \ \ $\rho(\overline W)\ =\ \overline V$, 
  \medskip\noindent
  (2) \ \  $\rho$ is an embedding on a tubular neighborhood of $\partial \overline W$ in $S$
  and 
  \medskip\noindent
  (3) \ \  $\rho(\partial \overline W)$ is a union of components of $\g$.
  }
 
\medskip
\noindent
{\bf Proof of Theorem \CC.1.}
 We assume $n =2$.  The case of general $n$ is similar.

 We first note that $V$ has finite area and finitely many irreducible components
 $V_1,...,V_\ell$.  This follows from work of Shiffman, but can be seen directly
 as follows.  Choose any $p\in \bbp^2-\overline V$ and let $\pi : \bbp^2-\{p\}\to \bbp^1$
 be projection.  Then $\pi\bigr|_V$ is finitely sheeted over $\bbp^1-\pi(\g)$, and therefore
 $V$ has finitely many components.
 In fact $\pi\bigr|_V$  must also be finitely sheeted over all of $\bbp^1$. To see this note
 that $V$ can contain no fibre of   $\pi$ since $p\notin \overline V =V\cup\g$. Hence, the intersection
 $\pi^{-1}(x)\cap V$ for $x\in \pi(\g)$ is at most countable.  If it were infinite, one easily sees that the
 sheeting  number in contiguous domains of $\bbp^1-\pi(\g)$ is unbounded. Choosing two distinct
 such projections and an easy estimate shows that the integral of the projective  K\"ahler form
 on $V$ is finite.
 
 Now each irreducible component $V_j$ defines a current $[V_j]$  by integration whose boundary is 
 supported in $\g$.  By the Federer Flat Support Theorem [3, 4.1.15],
  $$
  \partial [V_j]=n_j[\g_j]
  $$
   where $\g_j \equiv \supp \, \partial [V_j]$  is a union of connected 
   components of $\g$  (appropriately oriented) and
   $n_j\geq0$ is a locally constant   integer-valued function on $\g_j$.  
   Order the $V_j$ so that $n_j\geq 1$ for $j=1,...,m$ and $n_j=0$
 (that is, $\partial [V_j]=0$) for   $j>m$.
 
 Since $\g$ is a regularly embedded real analytic curve, it has a complexification 
 $\Sigma\supset \g$ which is a union of regularly embedded closed complex analytic annuli.
 Let $\Sigma_j$ denote that part of $\Sigma$ which is  the complexification of 
 $\g_j$ for
 $j\leq m$. Write $\Sigma_j = \Sigma_j^+\cup \g_j\cup \Sigma_j^-$ where $\Sigma_j^\pm$ are
 the components of $\Sigma_j-\g_j$ with signs chosen so that $\Sigma^+$ is the ``outer strip'',
 that is, so that 
 $$
 \partial \Sigma_j^+ \ =\ \g_j^+  -  \g_j.
 $$
 Consider the current $[V_j^*] \equiv[V_j]+ n_j[\Sigma_j^+]$ which has 
 $$
 \partial [V_j^*]  = n_j [\partial \Sigma_j^+].
 $$ 
 The structure theorem of King [8] implies that $\supp [V_j^*] $
 is a 1-dimensional subvariety of $\bbp^2-  \g_j^+$. It follows that $V_j^*$
 is an analytic continuation of  $V_j$ and in particular
$$
n_j\equiv 1 \and 
 \Sigma_j^+\ss V_j.
 $$
 
Defining $\rho_j:S_j \to V_j^*$ to be the normalization of $V_j^*$ and setting $\overline W_j
=\rho^{-1}(\overline V_j)$ completes part (1).  

The remaining components of $V$ are algebraic curves.  If one of them, say $V_k$,
 contains a union $\g_k$  of  components of $\g$, then it contains the complexification
 of $\g_k$ which is a union regularly embedded of complex annuli. This proves part (2).
 Part (3) is obvious.\qed

\vfill\eject
\centerline{\bf References}

 \bigskip

  \noindent 
  \item{[1]}  H. Alexander and J. Wermer,  { Several Complex Variables and Banach Algebras},    
Springer Verlag, New York, 1998.

 \medskip
  \noindent 
  \item{[2]}  E. Bishop, {\sl Holomorphic completions, analytic continuations   and the interpolation of semi-norms}, Ann. of Math. {\bf 78} (1963), 468-500.

 \medskip    
\noindent
\item{[3]}  H. Federer, Geometric Measure  Theory,
 Springer--Verlag, New York, 1969.

 \medskip
  \noindent 
  \item{[4]}
 F. R. Harvey and H. B. Lawson, Jr., {\sl Projective hulls and the projective Gelfand transformation},  Asian J. Math. {\bf 10}, no. 3 (2006), 279-319.

  \medskip
  \noindent 
  \item{[5]}
 F. R. Harvey and H. B. Lawson, Jr., {\sl Projective linking and boundaries of 
 positive holomorphic chains in projective manifolds, Part I},  ArXiv:math.CV/0512379

 \medskip 
  \noindent 
  \item{[6]}
   F. R. Harvey and H. B. Lawson, Jr., {\sl Projective linking and boundaries of 
 positive holomorphic chains in projective manifolds, Part II},  ArXiv:math.CV/0608029
   
 \medskip
 
  \noindent 
  \item{[7]}
 F. R. Harvey and H. B. Lawson, Jr., {\sl  Relative holomorphic cycles and duality}, Stony Brook preprint, 2006.   
 \medskip

      \noindent 
  \item{[8]} J. King, {\sl The currents defined by analytic varieties}, Acta Math. {\bf 127}
  (1971), 185-220.

\medskip

  \noindent 
  \item{[9]}   J. Wermer,  {\sl    The hull of a curve in $\bbc^n$},    
Ann. of Math., {\bf  68}  (1958), 550-561.

 \medskip
   \noindent 
  \item{[10]} J. Wermer, { Banach Algebras and Several Complex Variables}, 
  Markham Publishing Co., 1971.

  \end
  
  (Intuitively, 
  $\PH \g \cap (D\times \D_R) $ only leaks out of the side  of the tin can.)

 $\PH \g \cap (D\times \D_R)$ does not
  meet   that part of the boundary corresponding to $D\times\partial \D_R$.   Otherwise said,